\documentclass[12pt]{amsart}
\usepackage[a4paper,margin=2.5cm,headsep=8mm]{geometry}
\usepackage{amsmath,amssymb,amsthm,mathtools}
\usepackage{booktabs}
\usepackage{array}
\usepackage{xcolor}
\usepackage{enumitem}
\usepackage[colorlinks=true,linkcolor=blue!60!black,citecolor=blue!60!black,urlcolor=blue!60!black]{hyperref}
\usepackage{orcidlink}
\usepackage{xurl}
\theoremstyle{plain}
\newtheorem{theorem}{Theorem}[section]
\newtheorem{lemma}[theorem]{Lemma}
\newtheorem{proposition}[theorem]{Proposition}

\newtheorem{conjecture}[theorem]{Conjecture}
\theoremstyle{definition}
\newtheorem{definition}[theorem]{Definition}
\newtheorem{example}[theorem]{Example}
\newtheorem{remark}[theorem]{Remark}
\newtheorem{problem}[theorem]{Problem}

\newcommand{\dd}{\mathrm{dd}}
\newcommand{\DD}{\mathrm{DD}}
\newcommand{\ddD}{\mathrm{dd}_{D}}
\newcommand{\DDD}{\mathrm{DD}_{D}}
\newcommand{\Bp}{\mathrm{Bp}}
\newcommand{\diam}{\operatorname{diam}}
\newcommand{\cart}{\mathbin{\square}}

\begin{document}

\title{Detour different-distance sets in graphs}

\author[A. Alsinai]{Ammar Alsinai\,\orcidlink{0000-0002-5221-0574}}
\address{Department of Mathematics and Computer Science, College of Engineering and Information Technology, Onaizah Colleges, Qassim, Saudi Arabia}
\email{aliiammar1985@gmail.com}
\thanks{Corresponding author: Ammar Alsinai (aliiammar1985@gmail.com).}

\author[H. Ahmed]{Hanan Ahmed\,\orcidlink{0000-0002-4008-4873}}
\address{Department of Mathematics, Ibb University, Ibb, Yemen}
\email{hananahmed1a@gmail.com}


\begin{abstract}
The detour distance $D(u,v)$ between two vertices of a connected graph is the length of a longest $u,v$-path. A set $S$ of vertices is a detour different-distance set if no vertex outside $S$ has the same detour distance to two vertices of $S$. We study the lower and upper detour different-distance numbers $\ddD(G)$ and $\DDD(G)$, the smallest and largest orders of maximal such sets; this is the detour analogue of the different-distance sets of Hedetniemi, Hedetniemi, Laskar and Mulder. The two notions coincide on trees, cycles, complete bipartite graphs and the Petersen graph, but differ in general: for instance, detour different-distance sets may contain triangles or induce $4$-cycles. We determine both numbers for many families, among them wheels, fans, windmills, hypercubes, helms, flowers, books, ladders, prisms, M\"obius ladders, crowns, gears, Dutch windmills and complete multipartite graphs, and for joins, coronas and lexicographic products. We prove that $\DDD(G)$ is at most the detour diameter; that $\DDD(G)=n-1$ if and only if $G$ is a path; that $\ddD(G)=\DDD(G)=1$ whenever two vertices are joined by Hamiltonian paths to all other vertices; that $\ddD(G)=\DDD(G)=2$ for Hamiltonian-laceable graphs and, more generally, for bipartite graphs in which each part contains two vertices whose detour distances to the other part are constant; and that $\DDD(G)\le\lfloor n/2\rfloor$ for Hamiltonian graphs.
\end{abstract}

\maketitle

\section{Introduction}\label{sec:intro}

All graphs are finite, simple and connected, of order $n\ge2$. We write $d(u,v)$ for the (geodesic) distance and use \cite{Harary,CL} for terminology.

The \emph{detour distance} $D(u,v)=D_G(u,v)$ between vertices $u$ and $v$ is the length of a longest $u,v$-path in $G$, with $D(u,u)=0$. It was introduced by Chartrand, Johns and Tian \cite{CJT}, who showed that it is a metric on $V(G)$; it has since been studied extensively (see \cite{CL} and the references there). The \emph{detour eccentricity} of $v$ is $e_D(v)=\max_uD(u,v)$, and the \emph{detour diameter} is $\diam_D(G)=\max_{u,v}D(u,v)$. Clearly $d(u,v)\le D(u,v)\le n-1$.

Hedetniemi, Hedetniemi, Laskar and Mulder \cite{HHLM} introduced \emph{different-distance sets}: a set $S$ with $\emptyset\ne S\ne V(G)$ such that no vertex outside $S$ has the same distance to two vertices of $S$. The lower and upper different-distance numbers $\dd(G)$ and $\DD(G)$ are the minimum and maximum orders of maximal such sets. In this paper we replace the geodesic distance by the detour distance.

\begin{definition}
A set $S\subseteq V(G)$ with $\emptyset\ne S\ne V(G)$ is a \emph{detour different-distance set} (DDD set) if $D(w,u)\ne D(w,v)$ for every $w\notin S$ and all distinct $u,v\in S$. The \emph{lower} and \emph{upper detour different-distance numbers} $\ddD(G)$ and $\DDD(G)$ are the minimum and the maximum order of a maximal DDD set. A two-element DDD set is a \emph{DDD pair}.
\end{definition}

Every singleton is a DDD set, so $1\le\ddD(G)\le\DDD(G)\le n-1$.

Throughout, as in \cite{HHLM}, \emph{maximal} is meant with respect to inclusion: a (detour) different-distance set $S$ is maximal if no (detour) different-distance set properly contains $S$. Since these sets are not closed under taking subsets, this is stronger than requiring only that $S\cup\{v\}$ fails to be such a set for every $v\notin S$. For instance, in the path $v_1\cdots v_5$ the set $\{v_1,v_4\}$ is a different-distance set to which no single vertex can be added, but it is contained in the different-distance set $\{v_1,v_2,v_3,v_4\}$; similarly, in the book $B_n$ of Theorem~\ref{thm:book} the spine $\{a,b\}$ is a DDD pair to which no single vertex can be added, but it lies in every page. All values in this paper, and those quoted in Theorems~\ref{pre:HHLM} and~\ref{pre:cycles}, refer to maximality with respect to inclusion.

\begin{example}\label{ex:paw}
Let $G$ be the graph with vertices $v_1,\dots,v_4$ and edges $v_1v_2,v_2v_3,v_3v_4,v_4v_2$ (a triangle with a pendant vertex). Its detour distances are $D(v_1,v_2)=1$, $D(v_1,v_3)=D(v_1,v_4)=3$ and $D=2$ within the triangle. The maximal DDD sets are $\{v_1,v_2\}$, $\{v_1,v_3\}$ and $\{v_1,v_4\}$, so $\ddD(G)=\DDD(G)=2$. In the geodesic setting the maximal different-distance sets of this graph are $\{v_1,v_2\}$, $\{v_3\}$ and $\{v_4\}$ (a different-distance set contains no two vertices of a triangle \cite{HHLM}, and $v_1$ is equidistant from $v_3$ and $v_4$), so $\dd(G)=1$.
\end{example}

\subsection{Known results}
We now recall, with their sources, the known results used in this paper; later we refer to them by their numbers here.

\begin{theorem}[\cite{CJT}]\label{pre:metric}
The detour distance $D$ is a metric on $V(G)$.
\end{theorem}

A \emph{pendant path} is a path with exactly one end of degree $1$ in $G$ and all internal vertices of degree $2$ in $G$; a \emph{bridging path} is a path whose internal vertices have degree $2$ in $G$, all of whose edges are bridges, and whose deletion disconnects $G$. $\Bp(G)$ is the order of a longest pendant or bridging path.

\begin{theorem}[\cite{HHLM}]\label{pre:HHLM}
\begin{enumerate}[label=\textup{(\roman*)}]
\item For $n\ge3$, $\dd(P_n)=n-1$ if $n-1$ is a power of $2$, and $\dd(P_n)=\frac{n-1}{r}+1$ otherwise, where $r$ is the largest odd prime divisor of $n-1$; and $\DD(P_n)=n-1$.
\item $\dd(K_{m,r})=\DD(K_{m,r})=2$ for $m+r\ge3$, and the maximal different-distance sets of $K_{m,r}$ are its edges.
\item $\dd=\DD=1$ for the Petersen graph.
\end{enumerate}
\end{theorem}

\begin{theorem}[\cite{HHLM}]\label{pre:cycles}
For $n\ge3$: $\dd(C_n)=\DD(C_n)=1$ if $n$ is an odd prime; for odd composite $n$, $\dd(C_n)=n/P$ and $\DD(C_n)=n/r$, where $P$ and $r$ are the largest and smallest prime divisors of $n$; for even $n$, $\dd(C_n)=2$ and $\DD(C_n)=n/r$, where $r$ is the smallest odd prime divisor of $n$ ($\DD(C_n)=2$ if $n$ is a power of $2$).
\end{theorem}

In \cite{HHLM} the value of $\dd(C_n)$ for odd composite $n$ is written as $n/t$ with $t$ a largest odd divisor of $n$; as the proof there shows (the maximal different-distance sets of $C_n$ are the sets of $n/p$ equally spaced vertices, $p$ a prime divisor of $n$), $t$ is to be read as the largest prime divisor $P$. For example, $\{v_0,v_3,v_6\}$ is a different-distance set of $C_9$ containing $v_0$; so no singleton is maximal, and $\dd(C_9)=\DD(C_9)=3$.

A balanced bipartite graph with parts $A,B$ is \emph{Hamiltonian-laceable} if any two vertices from different parts are joined by a Hamiltonian path.

\begin{theorem}[\cite{Havel}]\label{pre:havel}
The hypercube $Q_k$ is Hamiltonian-laceable for every $k\ge1$.
\end{theorem}

\begin{theorem}[\cite{Ore}]\label{pre:ore}
If $\deg u+\deg v\ge n+1$ for every two non-adjacent vertices $u,v$ of a graph $G$ of order $n$, then $G$ is Hamiltonian-connected.
\end{theorem}

\subsection*{Outline}
In Section~\ref{sec:families} we determine $\ddD$ and $\DDD$ for many families: complete graphs, wheels, fans, friendship and windmill graphs, complete bipartite graphs, hypercubes, the Petersen graph, helms, flowers, books, ladders, prisms, M\"obius ladders, crowns, gears, Dutch windmills, complete multipartite graphs, paths and cycles, and (for $\DDD$) lollipops and barbells, as well as for joins, coronas and lexicographic products. These examples show the new phenomena of the detour setting: the triangle and opposite-pair lemmas of \cite{HHLM} fail, and maximal DDD sets may induce triangles or $4$-cycles. In Section~\ref{sec:bounds} we prove the general results used there: a counting bound in terms of detour distances, the characterization of $\DDD(G)=n-1$, the behaviour of pendant and bridging paths, criteria via Hamiltonian paths, a criterion for bipartite graphs, and a comparison principle between the geodesic and detour settings.

\section{Detour different-distance numbers of some families of graphs}\label{sec:families}

In this section we use the general results of Section~\ref{sec:bounds}.

\subsection{Graphs in which the two settings agree}

\begin{theorem}\label{thm:agree}
For trees, cycles, complete bipartite graphs and the Petersen graph, the DDD sets are exactly the different-distance sets. Consequently:
\begin{enumerate}[label=\textup{(\roman*)}]
\item $\ddD(T)=\dd(T)$ and $\DDD(T)=\DD(T)$ for every tree $T$; for paths these are given by Theorem~\ref{pre:HHLM}(i);
\item $\ddD(C_n)=\dd(C_n)$ and $\DDD(C_n)=\DD(C_n)$, given by Theorem~\ref{pre:cycles};
\item $\ddD(K_{m,r})=\DDD(K_{m,r})=2$ for $m+r\ge3$;
\item $\ddD=\DDD=1$ for the Petersen graph.
\end{enumerate}
\end{theorem}

\begin{proof}
By Lemma~\ref{lem:compare} it suffices to show that for every vertex $w$, the partitions of $V\setminus\{w\}$ according to $d(w,\cdot)$ and according to $D(w,\cdot)$ coincide. In a tree, $D=d$. In $C_n$, $D(u,v)=n-d(u,v)$. In $K_{m,r}$ with parts $A,B$, $|A|=m\le r=|B|$, a longest path alternates between the parts, which gives $D(u,v)=2m-1$ for $u\in A$, $v\in B$, $D(u,v)=2m-2$ for $u,v\in A$, and $D(u,v)=2m$ (if $m<r$) or $2m-2$ (if $m=r$) for $u,v\in B$; so for each $w$ the vertices of its own part and of the other part receive two different values, as for $d$. In the Petersen graph, adjacent vertices are joined by no Hamiltonian path (such a path together with the edge joining its ends would be a Hamiltonian cycle, and the Petersen graph is not Hamiltonian) but by a path of length $8$, and non-adjacent vertices by a Hamiltonian path; as the graph is distance-transitive, a check of one adjacent and one non-adjacent pair suffices. Hence $D=8$ at distance $1$ and $D=9$ at distance $2$. The values then follow from Theorems~\ref{pre:HHLM} and~\ref{pre:cycles}.
\end{proof}

\subsection{Graphs with $\ddD=\DDD=1$}

\begin{theorem}\label{thm:one}
$\ddD(G)=\DDD(G)=1$ for: the complete graphs $K_n$ ($n\ge2$); the wheels $W_n=K_1+C_{n-1}$ ($n\ge4$); the fans $F_n=K_1+P_{n-1}$ ($n\ge3$); the windmills $K_m^{(k)}$ ($k\ge2$ copies of $K_m$, $m\ge3$, sharing a vertex), in particular the friendship graphs.
\end{theorem}

\begin{proof}
We use Theorem~\ref{thm:universal}: it suffices to find two vertices from which Hamiltonian paths lead to every other vertex. In $K_n$ every vertex has this property.

In $W_n$ (where, by symmetry of the rim, every rim vertex behaves like $c_1$, so that $W_n$ is in fact Hamiltonian-connected) let $h$ be the hub and $c_1,\dots,c_{n-1}$ the rim in cyclic order. From $h$ to $c_i$: $h\,c_{i+1}c_{i+2}\cdots c_{i-1}c_i$. From $c_1$ to $h$: $c_1c_2\cdots c_{n-1}h$; from $c_1$ to $c_j$: $c_1\cdots c_{j-1}\,h\,c_{n-1}c_{n-2}\cdots c_j$. So $h$ and $c_1$ qualify.

In $F_n$ ($n\ge4$) let $h$ be the hub and $p_1\cdots p_{n-1}$ the path. From $p_1$ to $h$: $p_1p_2\cdots p_{n-1}h$; from $p_1$ to $p_j$: $p_1\cdots p_{j-1}\,h\,p_{n-1}p_{n-2}\cdots p_j$. So $p_1$ and, symmetrically, $p_{n-1}$ qualify. $F_3=K_3$.

In $K_m^{(k)}$ with centre $c$, $D(u,v)=m-1$ if $u,v$ lie in a common copy of $K_m$ and $D(u,v)=2m-2$ otherwise (Lemma~\ref{lem:blocks}). Let $S$ be a DDD set with $|S|\ge2$. If $c\notin S$, $c$ is at detour distance $m-1$ from all other vertices; so $c\in S$. Let $x\in S\setminus\{c\}$ lie in the copy $Q$. Every other non-central vertex $z$ of $Q$ satisfies $D(z,c)=D(z,x)=m-1$, so $Q\subseteq S$; then a vertex $y$ of another copy satisfies $D(y,x)=D(y,x')=2m-2$ for two non-central vertices $x,x'$ of $Q$, so $y\in S$. Hence $S=V$, a contradiction.
\end{proof}

\subsection{Graphs with $\ddD=\DDD=2$}

\begin{theorem}\label{thm:two}
\begin{enumerate}[label=\textup{(\roman*)}]
\item For every Hamiltonian-laceable graph $G$ with parts $A,B$ and $n\ge4$, the maximal DDD sets are the pairs $\{a,b\}$ with $a\in A$, $b\in B$; so $\ddD(G)=\DDD(G)=2$, and there are $|A||B|$ maximal DDD sets. In particular this holds for the hypercubes $Q_k$, $k\ge2$, and for $K_{m,m}$, $m\ge2$.
\item For the helm $H_n$ ($n\ge3$), obtained from $W_{n+1}$ by attaching a pendant vertex $\ell_i$ to each rim vertex $c_i$, the maximal DDD sets are the $n(n+1)$ pairs consisting of a pendant vertex and a vertex of the wheel. So $\ddD(H_n)=\DDD(H_n)=2$.
\end{enumerate}
\end{theorem}

\begin{proof}
(i) is Theorem~\ref{thm:laceable} with Theorem~\ref{pre:havel}; $K_{m,m}$ is clearly Hamiltonian-laceable.

(ii) Let $W$ be the wheel with hub $h$. Since $W$ is a block of order $n+1$ in which every vertex lies on Hamiltonian paths to all others (Theorem~\ref{thm:one}), Lemma~\ref{lem:blocks} gives the detour distances: $D=n$ between two vertices of $W$; $D(c_i,\ell_i)=1$; $D(\ell_i,x)=n+1$ for $x\in V(W)\setminus\{c_i\}$; and $D(\ell_i,\ell_j)=n+2$. Hence the values seen from a vertex $z$ are: from $h$, $n$ on the rim and $n+1$ on the pendant vertices; from $c_i$, $n$ on $V(W)\setminus\{c_i\}$, $1$ on $\ell_i$ and $n+1$ on the other pendant vertices; from $\ell_i$, $1$ on $c_i$, $n+1$ on the rest of $W$ and $n+2$ on the other pendant vertices.

Let $S$ be a DDD set. If $S$ contains two vertices of $W$, then every vertex of $W$ outside $S$ sees them at the same value $n$, so $V(W)\subseteq S$, and then every pendant vertex $\ell_k$ sees two vertices of $V(W)\setminus\{c_k\}$ at the same value $n+1$, so $S=V$. If $S$ contains two pendant vertices $\ell_i,\ell_k$, then $h$ and every $c_j$ with $j\ne i,k$ see them at the same value, so $S$ contains two vertices of $W$. Hence $|S|\le2$, with at most one vertex of $W$ and at most one pendant vertex. Conversely, for $\ell_i$ and $x\in V(W)$, the list of values above shows that no vertex outside $\{\ell_i,x\}$ sees them at the same value. So the maximal DDD sets are these $n(n+1)$ pairs, and every vertex lies in one of them.
\end{proof}

\subsection{Flowers and books: DDD sets inducing triangles and $4$-cycles}

The flower $Fl_n$ is obtained from the helm $H_n$ by joining every pendant vertex $\ell_i$ also to the hub $h$. The book $B_n=K_{1,n}\cart K_2$ consists of a spine edge $ab$ and $n$ pages, the $4$-cycles $a\,x_i\,y_i\,b\,a$.

\begin{lemma}\label{lem:flowerD}
Let $n\ge4$. In $Fl_n$ the detour distances are $D(h,c_i)=n+1$, $D(h,\ell_i)=n+2$, $D(c_i,c_j)=n+2$, $D(c_i,\ell_i)=n+2$, $D(c_i,\ell_j)=n+3$ and $D(\ell_i,\ell_j)=n+4$ ($i\ne j$).
\end{lemma}

\begin{proof}
Upper bounds. A pendant vertex $\ell_k$ has the two neighbours $c_k,h$; if it is an internal vertex of a path, both are its neighbours on the path. As $h$ has at most two neighbours on a path, a path contains at most two internal pendant vertices, at most one if $h$ is an end, and none if $h$ is not on it. Since $Fl_n$ has $n+1$ vertices other than pendant ones, a $u,v$-path has at most $n+1$ such vertices, plus at most the number of internal pendant vertices allowed, plus the ends that are pendant. This gives at most $n+2$ vertices for $\{u,v\}=\{h,c_i\}$, $n+3$ for $\{h,\ell_i\}$ and for $\{c_i,c_j\}$, $n+4$ for $\{c_i,\ell_j\}$ and $n+5$ for $\{\ell_i,\ell_j\}$. For $\{c_i,\ell_i\}$, the end $\ell_i$ is preceded either by $c_i$ (then the path is the edge $c_i\ell_i$) or by $h$; in the latter case $h$ is internal with $\ell_i$ as one path neighbour, so at most one further pendant vertex is internal, and the path has at most $n+3$ vertices.

Lower bounds (indices modulo $n$). The two paths $h\,\ell_{i+1}c_{i+1}c_{i+2}\cdots c_{i}$ and $h\,\ell_{i+1}c_{i+1}c_{i+2}\cdots c_{i}\,\ell_i$ have lengths $n+1$ and $n+2$, and the path $c_i\,c_{i+1}\cdots c_{i-1}\,\ell_{i-1}\,h\,\ell_i$ has length $n+2$. For $i\ne j$ let
\begin{align*}
R(i,j)&=c_i\,c_{i-1}\cdots c_{j+1}\,\ell_{j+1}\,h\,\ell_{i+1}\,c_{i+1}\,c_{i+2}\cdots c_j,\\
R'(i,j)&=c_i\,c_{i+1}\cdots c_{j-1}\,\ell_{j-1}\,h\,\ell_{i-1}\,c_{i-1}\,c_{i-2}\cdots c_j ;
\end{align*}
each passes through all rim vertices, $h$ and two pendant vertices, so it has length $n+2$. Appending $\ell_j$ to $R(i,j)$ (if $j\ne i+1$) or to $R'(i,j)$ (if $j=i+1$; then its internal pendant vertices $\ell_{i-1},\ell_i$ differ from $\ell_j$ since $n\ge3$) gives a $c_i,\ell_j$-path of length $n+3$. Adding $\ell_i$ at the start and $\ell_j$ at the end of $R(i,j)$ or $R'(i,j)$ gives an $\ell_i,\ell_j$-path of length $n+4$ unless $i$ and $j$ are cyclically consecutive. If $j=i+1$, take $k\notin\{i,i+1,i+2\}$ (possible since $n\ge4$) and the path
\[
\ell_i\,c_i\,c_{i-1}\cdots c_k\,\ell_k\,h\,\ell_{k-1}\,c_{k-1}\cdots c_{i+1}\,\ell_{i+1},
\]
which lists the rim backwards from $c_i$ to $c_{i+1}$ and has length $n+4$.
\end{proof}

\begin{theorem}\label{thm:flower}
For $n\ge4$ the maximal DDD sets of the flower $Fl_n$ are the $n$ triangles $\{h,c_i,\ell_i\}$, so $\ddD(Fl_n)=\DDD(Fl_n)=3$. For $n=3$, $\ddD(Fl_3)=1$ and $\DDD(Fl_3)=2$.
\end{theorem}

\begin{proof}
By Lemma~\ref{lem:flowerD}, the values seen from a vertex are: from $h$, $n+1$ on the rim and $n+2$ on the pendant vertices; from $c_j$, $n+1$ on $h$, $n+2$ on $\{c_k:k\ne j\}\cup\{\ell_j\}$ and $n+3$ on $\{\ell_k:k\ne j\}$; from $\ell_j$, $n+2$ on $\{h,c_j\}$, $n+3$ on $\{c_k:k\ne j\}$ and $n+4$ on $\{\ell_k:k\ne j\}$.

The triangle $T_i=\{h,c_i,\ell_i\}$ is a DDD set: $c_j$ ($j\ne i$) sees $h,c_i,\ell_i$ at $n+1,n+2,n+3$, and $\ell_j$ at $n+2,n+3,n+4$.

Let $S$ be a DDD set. If $c_i,c_k\in S$ ($i\ne k$), then $h$ and every $c_j\notin S$ see them at equal values, so $h$ and all rim vertices lie in $S$, and then every $\ell_j$ sees two rim vertices other than $c_j$ at $n+3$; hence $S=V$. If $\ell_i,\ell_k\in S$, then $h$ sees them at $n+2$ and every $c_j$, $j\ne i,k$, at $n+3$, so $S$ contains at least two rim vertices (as $n\ge4$), which was excluded. If $c_i,\ell_k\in S$ with $k\ne i$, then $c_k$ sees them at $n+2$, so $c_k\in S$, again two rim vertices. Finally, if $h,c_i\in S$, then $\ell_i$ sees them at $n+2$, so $\ell_i\in S$. Consequently every DDD set of order at least $2$ is contained in some $T_i$, and the $T_i$ are maximal. The case $n=3$ is a finite check.
\end{proof}

\begin{lemma}\label{lem:bookD}
Let $n\ge3$. In the book $B_n$ the detour distances are $D(a,b)=3$, $D(a,x_i)=D(b,y_i)=5$, $D(a,y_i)=D(b,x_i)=4$, $D(x_i,y_i)=5$, $D(x_i,x_j)=D(y_i,y_j)=6$ and $D(x_i,y_j)=7$ ($i\ne j$).
\end{lemma}

\begin{proof}
A path meets each page $\{x_k,y_k\}$ in a subpath. If a page contains no end of the path but meets it, the path enters the page from $a$ and leaves it to $b$ (or conversely), so it uses both $a$ and $b$; hence at most one page without an end of the path meets it. This bounds the number of vertices of a $u,v$-path by $2$ (for $a,b$) plus $2$ (one further page) plus the vertices of the pages containing $u$ and $v$, and, for $u=a$, by the fact that a path from $a$ into a page and on to $b$ cannot return. One obtains the upper bounds $3,5,4,5,6,7$, attained by $a\,x_i\,y_i\,b$; $a\,x_j\,y_j\,b\,y_i\,x_i$; $a\,x_j\,y_j\,b\,y_i$; $x_i\,a\,x_j\,y_j\,b\,y_i$; $x_i\,y_i\,b\,y_k\,x_k\,a\,x_j$; and $x_i\,y_i\,b\,y_k\,x_k\,a\,x_j\,y_j$ with $k\notin\{i,j\}$ (here $n\ge3$ is used).
\end{proof}

\begin{theorem}\label{thm:book}
For $n\ge3$ the maximal DDD sets of the book $B_n$ are the $n$ pages $\{a,b,x_i,y_i\}$, which induce $4$-cycles, and the $n(n-1)$ pairs $\{x_i,y_j\}$ with $i\ne j$. Hence $\ddD(B_n)=2$ and $\DDD(B_n)=4$. For $n=2$, $\ddD(B_2)=\DDD(B_2)=2$.
\end{theorem}

\begin{proof}
By Lemma~\ref{lem:bookD}, the values seen from a vertex are: from $a$, $3$ on $b$, $4$ on the $y$'s, $5$ on the $x$'s; from $x_i$, $4$ on $b$, $5$ on $\{a,y_i\}$, $6$ on the other $x$'s, $7$ on the other $y$'s; symmetrically for $b$ and $y_i$ (exchange $a\leftrightarrow b$, $x\leftrightarrow y$).

A page $P_i=\{a,b,x_i,y_i\}$ is a DDD set: $x_j$ ($j\ne i$) sees $b,a,x_i,y_i$ at $4,5,6,7$, and $y_j$ sees $a,b,y_i,x_i$ at $4,5,6,7$. A pair $\{x_i,y_j\}$, $i\ne j$, is a DDD set: $a$ sees it at $5,4$, $b$ at $4,5$, $x_j$ at $6,5$, $y_i$ at $5,6$, and $x_k$, $y_k$ ($k\ne i,j$) at $6,7$ and $7,6$.

Let $S$ be a DDD set with $|S|\ge2$. Suppose $x_i,x_j\in S$ with $i\ne j$, and let $k\notin\{i,j\}$. Then $a$ sees $x_i,x_j$ at $5$, $x_k$ at $6$ and $y_k$ at $7$, so $a,x_k,y_k\in S$; $b$ sees two $x$'s at $4$, so $b\in S$; and $y_i$ (resp.\ $y_j$) sees $x_j,x_k$ (resp.\ $x_i,x_k$) at $7$, so $y_i,y_j\in S$. Repeating with every $k$, $S=V$, a contradiction. So $S$ contains at most one $x$ and, symmetrically, at most one $y$; thus $S\subseteq\{a,b,x_i,y_j\}$ for some $i,j$. If $i\ne j$ and $x_i,y_j\in S$, then $a\notin S$ (otherwise $x_j$ sees $a$ and $y_j$ at $5$) and $b\notin S$ (otherwise $y_i$ sees $b$ and $x_i$ at $5$), so $S=\{x_i,y_j\}$. In all other cases $S$ is contained in a page. Hence the maximal DDD sets are the pages and the pairs $\{x_i,y_j\}$, $i\ne j$. For $n=2$ the claim is a finite check.
\end{proof}

\subsection{Lollipops and barbells}

\begin{theorem}\label{thm:lollipop}
Let $L(m,n)$ ($m\ge3$, $n\ge1$) be the lollipop ($K_m$ with a path $a\,p_1\cdots p_n$ attached at $a\in V(K_m)$) and $Bb(m,n)$ ($m\ge3$, $n\ge0$) the barbell (two copies of $K_m$ with vertices $a,b$ joined by a path $a\,p_1\cdots p_n\,b$). Then
\[
\DDD(L(m,n))=n+1,\qquad \DDD(Bb(m,n))=n+2.
\]
\end{theorem}

\begin{proof}
By Lemma~\ref{lem:blocks}, the detour distance between two vertices of a clique is $m-1$, and detour distances add along the path. Let $X=V(K_m)\setminus\{a\}$ in $L(m,n)$, and $X_A,X_B$ the analogous sets in $Bb(m,n)$; each has at least two vertices.

\emph{Claim: a DDD set $S$ contains at most one vertex of $X$} (of $X_A$, of $X_B$). If $x,x'\in S\cap X$, every vertex of $X\cup\{a\}$ outside $S$ sees them at $m-1$, so $X\cup\{a\}\subseteq S$; then $p_1$ sees $x,x'$ at $m$, so $p_1\in S$, and inductively every vertex lies in $S$.

Lollipop: the path $a\,p_1\cdots p_n$ is a pendant path, a DDD set of order $n+1$ (Lemma~\ref{lem:pendant}). Conversely, choose $w\in X\setminus S$, which exists by the Claim since $|X|=m-1\ge2$; the detour distances from $w$ to the other vertices take the $n+1$ values $m-1$ (on $X\cup\{a\}$) and $m-1+j$ (on $p_j$), so $|S|\le n+1$ by Lemma~\ref{lem:values}.

Barbell: the path $a\,p_1\cdots p_n\,b$ is a bridging path of order $n+2$. Choose $w\in X_A\setminus S$ (possible by the Claim, as $|X_A|\ge2$); the detour distances from $w$ take the $n+3$ values $m-1$ (on $X_A\cup\{a\}$), $m-1+j$ (on $p_j$), $m+n$ (on $b$) and $2m-2+n+1$ (on $X_B$). If $|S|=n+3$, then $S$ contains $b$ and a vertex $y\in X_B$; a vertex $w'\in X_B\setminus S$ sees $b$ and $y$ at $m-1$, a contradiction. So $|S|\le n+2$.
\end{proof}

\subsection{Ladders, prisms, M\"obius ladders, crowns and gears}

We write the vertices of the ladder $L_n=P_n\cart K_2$, the prism $C_n\cart K_2$ and the M\"obius ladder $M_n$ as $u_i,v_i$ (the $i$-th \emph{column}), $0\le i\le n-1$, with rungs $u_iv_i$ and rails $u_iu_{i+1}$, $v_iv_{i+1}$ ($0\le i\le n-2$). The prism has in addition the edges $u_{n-1}u_0$, $v_{n-1}v_0$, and $M_n$ the twisted edges $u_{n-1}v_0$, $v_{n-1}u_0$; thus $M_n$ is the cycle $u_0\cdots u_{n-1}v_0\cdots v_{n-1}u_0$ together with its $n$ long diagonals. The ladder and the even prisms are bipartite, and $M_n$ is bipartite exactly when $n$ is odd. The crown $S_n^0$ ($n\ge3$) is $K_{n,n}$ with parts $\{a_i\}$, $\{b_i\}$ minus the perfect matching $\{a_ib_i\}$. The gear $G_m$ ($m\ge3$) is the cycle $c_0c_1\cdots c_{2m-1}c_0$ together with a hub $h$ adjacent to $c_0,c_2,\dots,c_{2m-2}$.

\begin{lemma}\label{lem:zigzag}
In $C_n\cart K_2$ $(n\ge3)$ and in $M_n$ $(n\ge2)$ there are Hamiltonian paths from $u_0$ to $v_0$ and, for $1\le j\le n-1$,
\begin{enumerate}[label=\textup{(\roman*)}]
\item in $C_n\cart K_2$: to $u_j$ if $n-j$ is odd, and to $v_j$ if $n-j$ is even;
\item in $M_n$: to $v_j$ if $n-j$ is odd, and to $u_j$ if $n-j$ is even.
\end{enumerate}
\end{lemma}

\begin{proof}
The path $u_0u_1\cdots u_{n-1}v_{n-1}\cdots v_0$ ends at $v_0$. For $1\le j\le n-1$ start with the path
\[u_0u_1\cdots u_{j-1}\,v_{j-1}v_{j-2}\cdots v_0,\]
which covers columns $0,\dots,j-1$, and cross to column $n-1$, entering it at $v_{n-1}$ in the prism and at $u_{n-1}$ in $M_n$. Then traverse the columns $n-1,n-2,\dots,j$ in this order, each through its rung, moving from one column to the next along the rail of the row where we are (a \emph{zigzag}). Column $n-1-t$ is traversed in the same direction as column $n-1$ if $t$ is even and in the opposite direction if $t$ is odd; for column $j$, $t=n-1-j$, which gives the stated end vertices.
\end{proof}

\begin{theorem}\label{thm:ladders}
\begin{enumerate}[label=\textup{(\roman*)}]
\item For $n\ge2$, the maximal DDD sets of the ladder $L_n$ are the $n^2$ pairs of vertices from different parts; $\ddD(L_n)=\DDD(L_n)=2$.
\item For odd $n\ge3$, $\ddD(C_n\cart K_2)=\DDD(C_n\cart K_2)=1$; for even $n\ge4$, $C_n\cart K_2$ is Hamiltonian-laceable, its maximal DDD sets are the $n^2$ pairs from different parts, and $\ddD=\DDD=2$.
\item For even $n\ge2$, $\ddD(M_n)=\DDD(M_n)=1$; for odd $n\ge3$, $M_n$ is Hamiltonian-laceable, its maximal DDD sets are the $n^2$ pairs from different parts, and $\ddD=\DDD=2$.
\item For $n\ge3$, the crown $S_n^0$ is Hamiltonian-laceable; its maximal DDD sets are the $n^2$ pairs $\{a_i,b_j\}$ and $\ddD=\DDD=2$.
\item For $m\ge3$, the maximal DDD sets of the gear $G_m$ are the $m(m+1)$ pairs consisting of an even rim vertex $c_{2i}$ and a vertex of $\{h,c_1,c_3,\dots,c_{2m-1}\}$; $\ddD(G_m)=\DDD(G_m)=2$.
\end{enumerate}
\end{theorem}

\begin{proof}
(i) Colour $u_i$ by the parity of $i$ and $v_i$ by the opposite parity. We show that $u_0$ is saturating with constant $2n-1$. For $j\ge1$, traverse the columns $0,1,\dots,j-1$ by a zigzag starting at $u_0$ (so $u_0v_0v_1u_1u_2\cdots$); this enters column $j$ at $u_j$ if $j$ is even and at $v_j$ if $j$ is odd. Continue along that row to column $n-1$, cross the rung, and return along the other row to column $j$. The resulting Hamiltonian path ends at $v_j$ if $j$ is even and at $u_j$ if $j$ is odd, that is, at an arbitrary vertex of the other part; for $j=0$ use $u_0\cdots u_{n-1}v_{n-1}\cdots v_0$. By the symmetries of the ladder, the four corners $u_0,v_0,u_{n-1},v_{n-1}$ are saturating with constant $2n-1$, and two of them lie in each part. Apply Theorem~\ref{thm:saturating}.

(ii), (iii) The reflection $u_i\mapsto u_{-i}$, $v_i\mapsto v_{-i}$ (indices modulo $n$) is an automorphism of the prism fixing $u_0$. For odd $n$, Lemma~\ref{lem:zigzag}(i) and this reflection give Hamiltonian paths from $u_0$ to all other vertices, so $u_0$ is universal; since the prism is vertex-transitive, every vertex is universal and Theorem~\ref{thm:universal} applies. For even $n$, the vertices reached in Lemma~\ref{lem:zigzag}(i) are exactly those of the other part; by vertex-transitivity (automorphisms of a connected bipartite graph preserve the bipartition up to swapping the parts), the prism is Hamiltonian-laceable, and Theorem~\ref{thm:laceable} applies. For $M_n$ the map $x\mapsto -x$ on the cycle $u_0\cdots u_{n-1}v_0\cdots v_{n-1}$ (labelled $0,\dots,2n-1$) is an automorphism fixing $u_0$ that sends $u_j$ to $v_{n-j}$ and $v_j$ to $u_{n-j}$; the same argument with Lemma~\ref{lem:zigzag}(ii) applies, the parity conditions now being reversed.

(iv) Let $1\le j\le n-1$ and let $(p_0,\dots,p_{n-1})$ be a permutation of $\{0,\dots,n-1\}$ with $p_0=0$, $p_1=j$. Then
\[
a_{p_0}\,b_{p_2}\,a_{p_1}\,b_{p_3}\,a_{p_2}\,b_{p_4}\cdots a_{p_{n-1}}\,b_{p_{n+1}}\quad(\text{indices of }p\text{ modulo }n)
\]
is a Hamiltonian path from $a_0$ to $b_{p_1}=b_j$, since $p_{k+2}\ne p_k,p_{k+1}$ for $n\ge3$. By symmetry the crown is Hamiltonian-laceable; apply Theorem~\ref{thm:laceable}.

(v) The parts are $A=\{c_0,c_2,\dots,c_{2m-2}\}$ and $B=\{h,c_1,\dots,c_{2m-1}\}$, with $|A|=m<|B|$, so $D(x,y)\le2m-1$ for $x\in A$, $y\in B$. The paths
\begin{gather*}
h\,c_{2i+2}c_{2i+3}\cdots c_{2i-1}c_{2i},\qquad c_1c_0c_{2m-1}\cdots c_4\,h\,c_2,\\
c_1c_2\cdots c_{2i-2}\,h\,c_0c_{2m-1}c_{2m-2}\cdots c_{2i}\quad (2\le i\le m-1)
\end{gather*}
(indices modulo $2m$) have length $2m-1$ (each misses exactly one odd rim vertex) and show that $D(h,c_{2i})=2m-1$ for all $i$, and $D(c_1,c_{2i})=2m-1$ for $1\le i\le m-1$; the reflection $c_k\mapsto c_{2-k}$ fixes $c_1$ and $h$ and gives $D(c_1,c_0)=2m-1$. So $h$ and $c_1$ are saturating with constant $2m-1$. Rotations by an even number of steps are automorphisms; they show that every odd rim vertex is saturating, and also that $D(c_{2i},c_{2j+1})=D(c_{2i-2j},c_1)=2m-1$; together with $D(c_{2i},h)=2m-1$ this shows that the even rim vertices are saturating with constant $2m-1$ as well. Apply Theorem~\ref{thm:saturating}.
\end{proof}

\subsection{Dutch windmills and complete multipartite graphs}

The Dutch windmill $D_4^{(k)}$ consists of $k$ copies $c\,x_i\,y_i\,z_i\,c$ $(1\le i\le k)$ of $C_4$ sharing the vertex $c$.

\begin{theorem}\label{thm:dutch}
For $k\ge3$ the maximal DDD sets of $D_4^{(k)}$ are the $2k(k+1)$ pairs consisting of a vertex of $A=\{c,y_1,\dots,y_k\}$ and a vertex of $B=\{x_1,z_1,\dots,x_k,z_k\}$. Hence $\ddD(D_4^{(k)})=\DDD(D_4^{(k)})=2$.
\end{theorem}

\begin{proof}
In $C_4$ the detour distance is $3$ between adjacent and $2$ between opposite vertices. By Lemma~\ref{lem:blocks}, $D(c,x_i)=D(c,z_i)=3$, $D(c,y_i)=2$, $D(x_i,z_i)=2$, $D(x_i,y_i)=D(z_i,y_i)=3$, and $D(u,v)=D(u,c)+D(c,v)$ for $u,v$ in different copies. $A$ and $B$ are the parts of the bipartition. Let $S$ be a DDD set.

(a) If $c\in S$ and $y_i\in S$ for some $i$, then $x_i$ and $z_i$ see $c$ and $y_i$ at the same value $3$, so $x_i,z_i\in S$; then every vertex $w$ outside the $i$-th copy sees $x_i,z_i$ at $D(w,c)+3$, so $S=V$, a contradiction.

(b) If $S$ contains two vertices of $B$, then $c$ sees them at $3$, so $c\in S$; since $k\ge3$ there is a copy $j$ containing neither of them, and $y_j$ sees both at $2+3=5$, so $y_j\in S$, contradicting (a).

(c) If $S$ contains $y_i,y_j$ $(i\ne j)$, then $c$ sees them at $2$, so $c\in S$, contradicting (a); and $S$ cannot contain $c$ and some $y_i$ by (a).

So a DDD set has at most one vertex in each part, and Lemma~\ref{lem:parity} completes the proof.
\end{proof}

For the complete multipartite graph $K_{p_1,\dots,p_k}$ ($k\ge2$, $n=p_1+\dots+p_k\ge3$) let $P$ be a largest part, $p=|P|$, and $r=n-p$. For $k=2$ these are the complete bipartite graphs of Theorem~\ref{thm:agree}.

\begin{theorem}\label{thm:multipartite}
Let $G=K_{p_1,\dots,p_k}$ with $k\ge3$.
\begin{enumerate}[label=\textup{(\roman*)}]
\item If $p\ge r+1$, then $D(x,x')=2r$, $D(x,y)=2r-1$ and $D(y,y')=2r-2$ for $x,x'\in P$ and $y,y'\notin P$; the maximal DDD sets are the $pr$ pairs $\{x,y\}$ with $x\in P$, $y\notin P$, and $\ddD(G)=\DDD(G)=2$.
\item If $p\le r$, then $\ddD(G)=\DDD(G)=1$. In particular $\ddD=\DDD=1$ for the cocktail party graphs $K_{2,2,\dots,2}$ with at least three parts.
\end{enumerate}
\end{theorem}

\begin{proof}
(i) A path alternates between $P$ and $R=V\setminus P$ as long as possible: two vertices of $P$ on a path are separated by at least one vertex of $R$. A path between $x,x'\in P$ therefore contains at most $r+1$ vertices of $P$ and at most $r$ of $R$, and the path $x\,y_1\,x_1\,y_2\cdots y_r\,x'$ with $y_1,\dots,y_r$ the vertices of $R$ and $x_1,\dots,x_{r-1}\in P\setminus\{x,x'\}$ (available since $p\ge r+1$) attains $2r$. In the same way $D(x,y)=2r-1$ and $D(y,y')=2r-2$. Since these three values are distinct, for $x\in P$, $y\in R$ no vertex outside $\{x,y\}$ sees $x$ and $y$ at the same value, so $\{x,y\}$ is a DDD pair. If a DDD set $S$ contains $x,x'\in P$, every vertex of $R\setminus S$ sees them at $2r-1$, so $R\subseteq S$; since $k\ge3$, $r\ge2$, and every vertex of $P\setminus S$ sees two vertices of $R$ at $2r-1$, so $S=V$. If $S$ contains $y,y'\in R$, then every vertex of $P\setminus S$ sees them at $2r-1$, so $P\subseteq S$, which was just excluded. So the maximal DDD sets are the pairs $\{x,y\}$.

(ii) If $p<r$, two non-adjacent vertices lie in a common part of order $q\le p$, and their degree sum is $2(n-q)\ge2(n-p)\ge n+1$; so $G$ is Hamiltonian-connected by Theorem~\ref{pre:ore}, and Theorem~\ref{thm:universal} applies. If $p=r=n/2$, then $p\ge2$ (as $n\ge4$). Let $x\in P$ and $t\ne x$. A sequence of all vertices starting with $x$ is a Hamiltonian path if and only if consecutive vertices lie in different parts. If $t\in R$, alternate $x$, a vertex of $R$, a vertex of $P$, \dots, ending with $t$. If $t\in P$, write the $p$ vertices of $P$ with $x$ first and $t$ last and fill the $p-1$ gaps with the $r=p$ vertices of $R$, one gap receiving two vertices; since $R$ meets at least two parts ($k\ge3$), these two can be taken from different parts. So every vertex of $P$ is universal, and Theorem~\ref{thm:universal} applies.
\end{proof}

\subsection{Graph operations}

For graphs $G_1,G_2$ with $p_1$ and $p_2$ vertices, the join $G_1+G_2$ joins every vertex of $G_1$ to every vertex of $G_2$; the corona $G_1\circ G_2$ takes a copy $G_2^{x}$ of $G_2$ for every $x\in V(G_1)$ and joins $x$ to all vertices of $G_2^{x}$; and the lexicographic product $G[K_m]$ replaces every vertex $x$ of $G$ by a clique $Q_x$ of order $m$, all vertices of $Q_x$ and $Q_y$ being adjacent when $xy\in E(G)$. A graph is \emph{traceable} if it has a Hamiltonian path.

\begin{theorem}[Join]\label{thm:join}
Let $G_1,G_2$ be graphs (not necessarily connected) such that $G_2$ is traceable, $p_2\ge2$ and $1\le p_1\le p_2-1$. Then $\ddD(G_1+G_2)=\DDD(G_1+G_2)=1$. In particular this holds for $K_1+G$ with $G$ traceable of order at least $2$ (e.g.\ fans and wheels), for the complete split graphs $rK_1+K_s$ with $r<s$ (here $rK_1$ denotes $r$ isolated vertices), and for $rK_1+G$ with $G$ traceable of order larger than $r$. Neither hypothesis can be omitted: $K_{m,r}=mK_1+rK_1$ ($m,r\ge2$) and $K_1+(K_2\cup K_1)$ (Example~\ref{ex:paw}) have $\ddD=\DDD=2$ and fail traceability, while $4K_1+P_3$ (four independent vertices joined to a path of order $3$) and $K_{1,4}+K_2$ have $\ddD=\DDD=2$ and fail the size condition in both orders. The size condition is not necessary, however: $K_3+P_3$ fails it in both orders and has $\ddD=\DDD=1$.
\end{theorem}

\begin{proof}
Let $y_1y_2\cdots y_{p_2}$ be a Hamiltonian path of $G_2$. Every vertex of $G_1$ is adjacent to every $y_i$, so it can be inserted between any two consecutive vertices of a path of vertices of $G_2$. We show that $y_1$ is universal; by symmetry so is $y_{p_2}$, and Theorem~\ref{thm:universal} applies.

To end at $y_{p_2}$: insert the $p_1\le p_2-1$ vertices of $G_1$ into distinct gaps of $y_1\cdots y_{p_2}$. To end at $x\in V(G_1)$: insert the other $p_1-1$ vertices of $G_1$ into distinct gaps of $y_1\cdots y_{p_2}$ and append $x$. To end at $y_j$, $2\le j\le p_2-1$: take $y_1\cdots y_{j-1}\,x\,y_{p_2}y_{p_2-1}\cdots y_j$ for some $x\in V(G_1)$, and insert the remaining $p_1-1$ vertices of $G_1$ into the $p_2-2\ge p_1-1$ gaps between consecutive $y$'s of this path.
\end{proof}

\begin{theorem}[Corona]\label{thm:corona}
Let $G_1$ be Hamiltonian-connected of order $p_1\ge3$, and let $G_2$ have order $p_2$ such that $K_1+G_2$ is Hamiltonian-connected (e.g.\ $G_2=K_{p_2}$, or $G_2=C_{p_2}$ with $p_2\ge3$, or $G_2=W_{p_2}$).
\begin{enumerate}[label=\textup{(\roman*)}]
\item If $p_2=1$, the maximal DDD sets of $G_1\circ K_1$ are the $p_1^{2}$ pairs consisting of a pendant vertex and a vertex of $G_1$; so $\ddD=\DDD=2$.
\item If $p_2\ge2$ and $p_2\ne p_1-1$, the maximal DDD sets of $G_1\circ G_2$ are the $p_1p_2(p_1-1)$ pairs $\{w,y\}$ with $w\in V(G_2^{x})$ and $y\in V(G_1)\setminus\{x\}$; so $\ddD=\DDD=2$.
\item If $p_2=p_1-1$, then $\ddD(G_1\circ G_2)=\DDD(G_1\circ G_2)=1$.
\end{enumerate}
For example, $K_n\circ K_m$ ($n\ge3$) has $\ddD=\DDD=2$ if $m\ne n-1$, and $\ddD=\DDD=1$ if $m=n-1$.
\end{theorem}

\begin{proof}
A Hamiltonian-connected graph of order at least $3$ is $2$-connected, so the blocks of $G_1\circ G_2$ are $G_1$ and the cones $\{x\}\cup V(G_2^{x})$. By Lemma~\ref{lem:blocks}, with $X_x=V(G_2^{x})$:
$D(x,y)=p_1-1$ for $x\ne y$ in $G_1$; $D(w,x)=p_2$ and $D(w,w')=p_2$ for $w,w'\in X_x$; $D(w,y)=p_2+p_1-1$ for $w\in X_x$, $y\in V(G_1)\setminus\{x\}$; and $D(w,u)=2p_2+p_1-1$ for $w\in X_x$, $u\in X_y$, $y\ne x$.

Let $S$ be a DDD set. If $y,y'\in S\cap V(G_1)$, every vertex of $G_1$ outside $S$ sees them at $p_1-1$, so $V(G_1)\subseteq S$, and then every $w\in X_x\setminus S$ sees two vertices of $V(G_1)\setminus\{x\}$ at $p_2+p_1-1$; so $S=V$. Hence $|S\cap V(G_1)|\le1$. If $w,w'\in S$ lie in the same $X_x$ ($p_2\ge2$), then $x$ and the other vertices of $X_x$ see them at $p_2$, and every $y\in V(G_1)\setminus\{x\}$ at $p_2+p_1-1$, which puts two vertices of $G_1$ into $S$. If $w\in X_x$, $w'\in X_y$ with $x\ne y$, every $z\in V(G_1)\setminus\{x,y\}$ sees them at $p_2+p_1-1$; if $p_1\ge4$ this puts two vertices of $G_1$ into $S$, and if $p_1=3$ it puts the unique such $z$ into $S$, after which every vertex of $X_z$ sees $w,w'$ at $2p_2+p_1-1$ and $x$ sees $w'$ and a vertex of $X_z$ at $p_2+p_1-1$, again giving two vertices of $G_1$. So $S$ contains at most one vertex of $G_1$ and at most one vertex of $\bigcup_xX_x$, and $|S|\le2$.

(i) For $p_2=1$ write $X_x=\{\ell_x\}$, so $D(\ell_x,x)=1$, $D(\ell_x,y)=p_1$, $D(\ell_x,\ell_y)=p_1+1$. Consider $\ell_x$ and $y\in V(G_1)$. If $y=x$, a vertex $z\in V(G_1)\setminus\{x\}$ sees $\ell_x,x$ at $p_1,p_1-1$, and $\ell_z$ at $p_1+1,p_1$. If $y\ne x$, the vertex $x$ sees $\ell_x,y$ at $1,p_1-1$; $z\in V(G_1)\setminus\{x,y\}$ at $p_1,p_1-1$; $\ell_y$ at $p_1+1,1$; and $\ell_z$ ($z\ne x,y$) at $p_1+1,p_1$. As $p_1\ge3$, all these values differ, so all these pairs are DDD pairs.

(ii), (iii) For $p_2\ge2$, the pair $\{w,x\}$ with $w\in X_x$ is seen at $p_2$ by every other vertex of $X_x$, so it is not a DDD set. For $w\in X_x$ and $y\ne x$: the vertex $x$ sees them at $p_2$ and $p_1-1$; $z\in V(G_1)\setminus\{x,y\}$ at $p_2+p_1-1$ and $p_1-1$; $u\in X_x$ at $p_2$ and $p_2+p_1-1$; $u\in X_y$ at $2p_2+p_1-1$ and $p_2$; $u\in X_z$ at $2p_2+p_1-1$ and $p_2+p_1-1$. So $\{w,y\}$ is a DDD pair iff $p_2\ne p_1-1$. In case (ii) every vertex lies in such a pair; in case (iii) there are no DDD pairs at all, so every maximal DDD set is a singleton.
\end{proof}

\begin{theorem}[Lexicographic product]\label{thm:lex}
Let $G$ be a traceable graph of order $p\ge2$ and $m\ge2$. Then every vertex of the cliques corresponding to the two ends of a Hamiltonian path of $G$ is universal in $G[K_m]$; hence $\ddD(G[K_m])=\DDD(G[K_m])=1$. In particular this holds for $P_p[K_m]$, $C_p[K_m]$ and $K_p[K_m]$.
\end{theorem}

\begin{proof}
Let $x_1\cdots x_p$ be a Hamiltonian path of $G$, $Q_i=Q_{x_i}$, and $w\in Q_1$. It suffices to work in the spanning subgraph $P_p[K_m]$, in which consecutive cliques are completely joined. Let $t$ be a target. If $t\in Q_p$: traverse $Q_1$ (starting at $w$), $Q_2,\dots,Q_{p-1}$ completely, then $Q_p\setminus\{t\}$, then $t$. If $t\in Q_j$ with $1\le j<p$: traverse $Q_1,\dots,Q_{j-1}$ completely (starting at $w$; if $j=1$ start with $w$ only), then one vertex $a_j\in Q_j\setminus\{t,w\}$ (if $j>1$), then one vertex $a_i$ from each of $Q_{j+1},\dots,Q_p$; then return through the remaining vertices of $Q_p,Q_{p-1},\dots,Q_{j+1}$, which are non-empty since $m\ge2$, then through the remaining vertices of $Q_j$ other than $t$, and finish at $t$. All consecutive vertices lie in the same or in consecutive cliques, so this is a Hamiltonian $w,t$-path. Since $m\ge2$, there are at least two universal vertices, and Theorem~\ref{thm:universal} applies.
\end{proof}

\begin{remark}[Further operations]\label{rem:ops}
Subdivision: $S(T)$ is a tree for every tree $T$, and $S(C_n)=C_{2n}$, so their numbers follow from Theorem~\ref{thm:agree}. Cartesian products: $Q_k$ (Theorem~\ref{thm:two}), the books $K_{1,n}\cart K_2$ (Theorem~\ref{thm:book}), and the ladders and prisms $P_n\cart K_2$, $C_n\cart K_2$ (Theorem~\ref{thm:ladders}) were treated above.
\end{remark}

\begin{table}[h]
\centering
\small
\begin{tabular}{lcc}
\toprule
Graph & $\ddD$ & $\DDD$\\
\midrule
$K_n$, $W_n$, $F_n$, $K_m^{(k)}$ ($m\ge3$) & $1$ & $1$\\
Petersen graph & $1$ & $1$\\
$K_{m,r}$, $Q_k$ & $2$ & $2$\\
Helm $H_n$ & $2$ & $2$\\
Flower $Fl_n$ ($n\ge4$) & $3$ & $3$\\
Book $B_n$ ($n\ge3$) & $2$ & $4$\\
Lollipop $L(m,n)$ & -- & $n+1$\\
Barbell $Bb(m,n)$ & -- & $n+2$\\
Ladder $L_n$, crown $S_n^0$, gear $G_m$, $D_4^{(k)}$ ($k\ge3$) & $2$ & $2$\\
$C_n\cart K_2$, $n$ odd & $1$ & $1$\\
$C_n\cart K_2$, $n$ even; $M_n$, $n$ odd & $2$ & $2$\\
$M_n$, $n$ even & $1$ & $1$\\
$K_{p_1,\dots,p_k}$, $k\ge3$, $p\ge r+1$ & $2$ & $2$\\
$K_{p_1,\dots,p_k}$, $k\ge3$, $p\le r$ & $1$ & $1$\\
Trees, cycles & $\dd$ & $\DD$\\
\bottomrule
\end{tabular}
\caption{Detour different-distance numbers of some families (for trees and cycles they coincide with the geodesic numbers, see Theorem~\ref{thm:agree}).}
\label{tab:summary}
\end{table}

\section{Bounds and characterizations}\label{sec:bounds}

\subsection{A counting bound}

\begin{lemma}\label{lem:values}
Let $S$ be a DDD set and $w\notin S$. Then $D(w,\cdot)$ is injective on $S$; hence $|S|$ is at most the number of distinct values of $D(w,v)$, $v\ne w$. In particular $\DDD(G)\le e_D(w)$ for every vertex $w$ outside a maximum DDD set, and
\[
\DDD(G)\le\diam_D(G)\le n-1.
\]
\end{lemma}

\begin{proof}
This is the definition, together with $1\le D(w,v)\le e_D(w)$.
\end{proof}

\begin{theorem}\label{thm:path}
$\DDD(G)=n-1$ if and only if $G\cong P_n$.
\end{theorem}

\begin{proof}
For $P_n$ the detour distance equals the distance, and $\DD(P_n)=n-1$ (Theorem~\ref{pre:HHLM}(i)). Conversely, let $S$ be a DDD set of order $n-1$ and $V\setminus S=\{w\}$. By Lemma~\ref{lem:values}, $D(w,\cdot)$ maps $S$ bijectively onto $\{1,\dots,n-1\}$. Let $s$ satisfy $D(w,s)=n-1$ and let $w=x_0x_1\cdots x_{n-1}=s$ be a Hamiltonian path. Then $D(w,x_i)\ge i$, and since the values are a permutation of $1,\dots,n-1$, induction from $i=n-1$ downwards gives $D(w,x_i)=i$ for all $i$. If $x_ix_j\in E(G)$ with $j>i+1$, then $x_0\cdots x_ix_jx_{j-1}\cdots x_{i+1}$ is a $w,x_{i+1}$-path of length $j>i+1$; this includes the case $i=0$. Hence $G$ is the path $x_0\cdots x_{n-1}$.
\end{proof}

\subsection{Blocks, pendant and bridging paths}

\begin{lemma}\label{lem:blocks}
Let $u,v$ be vertices of $G$, and let $B_1,\dots,B_k$ be the blocks of a $u,v$-path in the block-cutvertex tree, with consecutive blocks meeting in the cut vertices $x_1,\dots,x_{k-1}$, $x_0=u$, $x_k=v$. Then
\[
D(u,v)=\sum_{i=1}^{k}D_{B_i}(x_{i-1},x_i).
\]
\end{lemma}

\begin{proof}
Every $u,v$-path passes through $x_1,\dots,x_{k-1}$ in this order and, between $x_{i-1}$ and $x_i$, stays in $B_i$, since leaving $B_i$ through a cut vertex would require returning through the same cut vertex. Conversely, longest paths in the blocks concatenate.
\end{proof}

\begin{lemma}\label{lem:pendant}
The vertex set of every pendant or bridging path is a DDD set, and, if $G\ne K_2$, the ends of every bridge form a DDD pair. Consequently $\Bp(G)\le\DDD(G)$.
\end{lemma}

\begin{proof}
Let $P$ be a pendant or bridging path. A vertex $w\notin P$ lies in a component of $G-V(P)$ attached to $P$ at an end $x$ of $P$, and by Lemma~\ref{lem:blocks}, $D(w,p)=D(w,x)+d_P(x,p)$ for $p\in V(P)$, since the edges of $P$ are bridges. These values are distinct. For a bridge $uv$, the same argument applies to the path $uv$.
\end{proof}

\begin{remark}
In contrast with the geodesic setting \cite{HHLM}, a DDD set may contain two vertices of a triangle (Theorem~\ref{thm:flower}), and a pair of opposite vertices of a $4$-cycle (Theorem~\ref{thm:book}); and a DDD set of order at least $3$ need neither be independent nor induce a path.
\end{remark}

\subsection{Hamiltonian paths}

\begin{theorem}\label{thm:universal}
Call $w$ \emph{universal} if $D(w,v)=n-1$ for every $v\ne w$, i.e.\ Hamiltonian paths join $w$ to all other vertices. If $G$ has two universal vertices, then $\ddD(G)=\DDD(G)=1$. In particular this holds for Hamiltonian-connected graphs.
\end{theorem}

\begin{proof}
A universal vertex outside a set $S$ sees all vertices of $S$ at the same value $n-1$; so every DDD set of order at least $2$ contains all universal vertices. If $u_1,u_2$ are universal and lie in $S$, every $w\notin S$ satisfies $D(w,u_1)=n-1=D(w,u_2)$. Hence $S=V$, a contradiction. So $G$ has no DDD set of order at least $2$; every singleton is then a maximal DDD set, and $\ddD(G)=\DDD(G)=1$.
\end{proof}

\begin{theorem}\label{thm:laceable}
Let $G$ be Hamiltonian-laceable with parts $A,B$ and $n\ge4$. Then the maximal DDD sets of $G$ are exactly the pairs $\{a,b\}$ with $a\in A$, $b\in B$, so $\ddD(G)=\DDD(G)=2$.
\end{theorem}

\begin{proof}
For $u,v$ in different parts $D(u,v)=n-1$; for $u,v$ in the same part, a path between them has one more vertex in their part than in the other, so $D(u,v)\le n-2$. If a DDD set $S$ contained two vertices of $A$, every vertex of $B\setminus S$ would see them at $n-1$, so $B\subseteq S$, and then every vertex of $A\setminus S$ would see two vertices of $B$ at $n-1$; so $S=V$. Hence a DDD set has at most one vertex in each part. Conversely, for $a\in A$, $b\in B$, a vertex $w\in A\setminus\{a\}$ sees $b$ at $n-1$ and $a$ at most $n-2$, and similarly for $w\in B$; so $\{a,b\}$ is a DDD pair.
\end{proof}

\subsection{Bipartite graphs}

\begin{lemma}\label{lem:parity}
Let $G$ be a connected bipartite graph with parts $A,B$ and $n\ge3$. For every vertex $w$ and all $a\in A$, $b\in B$ with $w\notin\{a,b\}$, the numbers $D(w,a)$ and $D(w,b)$ have different parities. Consequently every pair $\{a,b\}$ with $a\in A$, $b\in B$ is a DDD pair.
\end{lemma}

\begin{proof}
In a bipartite graph the length of every path between two vertices is even if they lie in the same part and odd otherwise. Hence $D(w,a)\not\equiv D(w,b)\pmod 2$, and no vertex outside $\{a,b\}$ sees $a$ and $b$ at the same value; since $n\ge3$, $\{a,b\}\ne V$.
\end{proof}

Let $G$ be bipartite with parts $A,B$. We call a vertex $w$ \emph{saturating} with constant $\lambda$ if $D(w,y)=\lambda$ for every vertex $y$ in the part not containing $w$. If $|A|\le|B|$, a path between the two parts has equally many vertices in both parts, so $D(w,y)\le2|A|-1$ for $w,y$ in different parts.

\begin{theorem}\label{thm:saturating}
Let $G$ be a connected bipartite graph with parts $A,B$, and suppose that $A$ contains two saturating vertices with a common constant and $B$ contains two saturating vertices with a common constant. Then the maximal DDD sets of $G$ are exactly the $|A||B|$ pairs $\{a,b\}$ with $a\in A$, $b\in B$. In particular $\ddD(G)=\DDD(G)=2$.
\end{theorem}

\begin{proof}
Let $s_A,s'_A\in A$ and $s_B,s'_B\in B$ be saturating, with $D(s_A,\cdot)=D(s'_A,\cdot)$ on $B$ and $D(s_B,\cdot)=D(s'_B,\cdot)$ on $A$; set $\lambda_A=D(s_A,y)$ $(y\in B)$ and $\lambda_B=D(s_B,x)$ $(x\in A)$. Let $S$ be a DDD set containing two vertices of $A$. If $s_B\notin S$, it sees these two vertices at the same value $\lambda_B$; so $s_B,s'_B\in S$. Now $S$ contains two vertices of $B$, and in the same way $s_A,s'_A\in S$. Let $w\notin S$. If $w\in A$, then $D(w,s_B)=\lambda_B=D(w,s'_B)$; if $w\in B$, then $D(w,s_A)=\lambda_A=D(w,s'_A)$. Both are impossible, so $S=V$, a contradiction. By symmetry no DDD set contains two vertices of $B$. Hence every DDD set has at most one vertex in each part, and by Lemma~\ref{lem:parity} all $|A||B|$ pairs $\{a,b\}$ are DDD pairs; they are maximal.
\end{proof}

Theorem~\ref{thm:laceable} is the special case in which every vertex is saturating with constant $n-1$. The hypothesis cannot be imposed on one part only: let $G$ be the $4$-cycle $u\,x\,v\,y\,u$ with a pendant vertex $\ell$ attached at $u$. The part $\{\ell,x,y\}$ contains the two saturating vertices $x,y$ (constant $3$), while in the part $\{u,v\}$ only $v$ is saturating; the maximal DDD sets are $\{\ell,u\}$, $\{u,x\}$, $\{u,y\}$, $\{\ell,x,v\}$ and $\{\ell,y,v\}$, so $\DDD(G)=3$.

\subsection{Hamiltonian graphs}

\begin{theorem}\label{thm:ham}
If $G$ is Hamiltonian, then $\DDD(G)\le\lfloor n/2\rfloor$.
\end{theorem}

\begin{proof}
Two vertices on a Hamiltonian cycle split it into two arcs, one of length at least $\lceil n/2\rceil$; so $D(w,v)\in\{\lceil n/2\rceil,\dots,n-1\}$ for $v\ne w$, a set of $\lfloor n/2\rfloor$ values. Apply Lemma~\ref{lem:values}.
\end{proof}

The bound is attained by $C_4$. The cycles $C_{3k}$, which have $\DDD=k$, suggest a stronger bound; see Conjecture~\ref{conj:ham}.

\subsection{Comparison with the geodesic setting}

\begin{lemma}\label{lem:compare}
Suppose that for every vertex $w$ and all $u,v\ne w$ we have $d(w,u)=d(w,v)$ if and only if $D(w,u)=D(w,v)$. Then the DDD sets of $G$ are exactly its different-distance sets, and $\ddD(G)=\dd(G)$, $\DDD(G)=\DD(G)$.
\end{lemma}

\begin{proof}
The two definitions then express the same condition.
\end{proof}

\begin{proposition}\label{prop:tree}
$D(u,v)=d(u,v)$ for all $u,v$ if and only if $G$ is a tree.
\end{proposition}

\begin{proof}
In a tree the $u,v$-path is unique. If $G$ has a cycle $C$, two adjacent vertices of $C$ satisfy $d=1$ and $D\ge|C|-1\ge2$.
\end{proof}

\begin{remark}\label{rem:compare}
In general neither number dominates the other. For the triangle with a pendant vertex, $\ddD=2>\dd=1$ (Example~\ref{ex:paw}), and for the complete multipartite graphs $K_{p_1,\dots,p_k}$ with $k\ge3$ and $p\ge r+1$ we have $\ddD=\DDD=2$ (Theorem~\ref{thm:multipartite}), while $\dd=\DD=1$ \cite{HHLM}. On the other hand, the prism $C_3\cart K_2=K_3\cart K_2$ has $\ddD=\DDD=1$ (Theorem~\ref{thm:ladders}), while $\dd=\DD=2$, since it is a Hamming graph with a coordinate of order $2$ \cite{HHLM}.
\end{remark}

\section{Concluding remarks and further problems}

The detour version of different-distance sets keeps several features of the geodesic one (pendant and bridging paths, the extremal role of paths, and complete agreement on trees and cycles), but loses the local structure: the triangle and opposite-pair lemmas fail and DDD sets need not be independent or induce paths. Hamiltonian properties take over the role played by distances.

\begin{conjecture}\label{conj:ham}
If $G$ is Hamiltonian of order $n\ge5$, then $\DDD(G)\le\max\{2,\lfloor n/3\rfloor\}$; equality holds for the cycles $C_{3k}$.
\end{conjecture}

The conjecture would improve Theorem~\ref{thm:ham} by a factor of about $2/3$: by the proof of that theorem, a maximal DDD set can use at most the $\lfloor n/2\rfloor$ values in $[\lceil n/2\rceil,n-1]$ seen from an outside vertex, and the conjecture asserts that only about two thirds of them can be used, as happens for the cycles.

\begin{problem}
Describe the subgraphs induced by DDD sets, in analogy with the description of the geodesic case in \cite{HHLM}.
\end{problem}

\begin{problem}
Determine $\ddD$ for lollipops, barbells and tadpoles, and both numbers for grids, sunlets, snakes and generalized Petersen graphs.
\end{problem}

\begin{problem}
Characterize the graphs with $\ddD(G)=\DDD(G)=1$, and those in which the DDD sets coincide with the different-distance sets.
\end{problem}

\section*{Funding}
This research received no external funding.

\section*{Data availability}
Data sharing is not applicable to this article, as no datasets were generated or analysed during the current study.

\section*{Author contributions}
Both authors contributed equally to this work.

\textbf{Ammar Alsinai:} Conceptualization, Methodology, Formal analysis, Investigation, Validation, Visualization, Writing -- original draft, Writing -- review \& editing.

\textbf{Hanan Ahmed:} Conceptualization, Methodology, Formal analysis, Investigation, Validation, Visualization, Writing -- original draft, Writing -- review \& editing.

All authors have read and approved the final version of the manuscript.

\section*{Conflict of interest}
The authors declare that they have no conflict of interest.

\textbf{Correspondence} should be addressed to Ammar Alsinai (A.A.).

\end{document}